\documentclass[11pt,a4paper]{article}
\usepackage{amsmath,amsthm}
\usepackage{url}
\usepackage[utf8]{inputenc}
\usepackage{float}
\usepackage{geometry}
\usepackage{graphicx}
\usepackage{amsfonts}
\usepackage{mathrsfs}
\usepackage[sort]{natbib}
\usepackage{fancyhdr}
\usepackage{verbatim}
\usepackage{enumerate}
\usepackage{paralist}
\usepackage[noend]{algpseudocode}
\usepackage{algorithmicx,algorithm}

\def\StudentName{Yijie Xu, Runqi Xu}			% insert here your name
\def\StudentNumber{27803035, 27803052}			% insert here your student number
\def\Proposer{Prof. Peter Sweby}		% insert here name of proposer

\newtheorem{theorem}{Theorem}[section]
\newtheorem*{defn}{Definition}

\title{Research on Interpolation and Data Fitting: Basis and Applications}
\author{\StudentName\\ \StudentNumber \\ \phantom{} \\ Department of Mathematics and Statistics \\ University of Reading }

\begin{document}
	\date{February 11, 2022}
	\maketitle

	%-----------------------------------------------------------------------------------ABSTRACT BEGINS
	
	\begin{abstract}
		\begin{normalsize}
		In the era of big data, we first need to manage the data, which requires us to find missing data or predict the trend, so we need operations including interpolation and data fitting. Interpolation is a process to discover deducing new data points in a range through known and discrete data points. When solving scientific and engineering problems, data points are usually obtained by sampling, experiments, and other methods. These data may represent a finite number of numerical functions in which the values of independent variables. According to these data, we often want to get a continuous function, i.e., curve; or more dense discrete equations are consistent with known data, while this process is called fitting.\\
		
		This article describes why the main idea come out logically and how to apply various method since the definitions are already written in the textbooks. At the same time, we give examples to help introduce the definitions or show the applications.\\
		
		We divide interpolation into several parts by their methods or functions for the structure. What comes first is the polynomial interpolation, which contains Lagrange interpolation and Newton interpolation, which are essential but critical. Then we introduce a typical stepwise linear interpolation —— Neville's algorithm. If we are concerned about the derivative, it comes to Hermite interpolation; if we focus on smoothness, it comes to cubic splines and Chebyshev nodes. Finally, in the Data fitting part, we introduce the most typical one: the Linear squares method, which needs to be completed by normal equations.
		\end{normalsize}
	\end{abstract}
	
	%-----------------------------------------------------------------------------------ABSTRACT ENDS

	\newpage
	
	\tableofcontents
	\listoffigures % NOT EDITED
	\listoftables % NOT EDITED
	
	\newpage
	\pagenumbering{arabic}
	%-----------------------------------------------------------------------------------REPORT BODY BEGINS

	\section{Interpolation by Polynomials}
	\subsection{Introduction}
	\large Put simply, \emph{polynomial interpolation} is the interpolation of a given data set through polynomials. Given a set of $n+1$ distinct data points $(x_i,y_i), i=0,1,\cdots,n$, we try to find a polynomial $p$, such that
	$$
	p(x_i)=y_i, i=0,1,\cdots,n.
	$$
	
	Note that the uniqueness theorem shows that such polynomial $p$ exists and is unique, which would be covered later.
	\subsection{Lagrange Interpolation}
	\subsubsection{Main Idea}
	Suppose we have three distinct points $(x_1,y_1), (x_2,y_2), (x_3,y_3)$, obviously we have to construct a $3$-order curve which could go across all three points. Rather than find it directly by its form as $f(x)=ax^2+bx+c, a\neq 0$, we want to follow Joseph-Louis Lagrange's thought: try to construct proper function by the composition of three functions related to these three points.
	\begin{figure}[h!] 
		\centering
		\includegraphics[width=14cm]{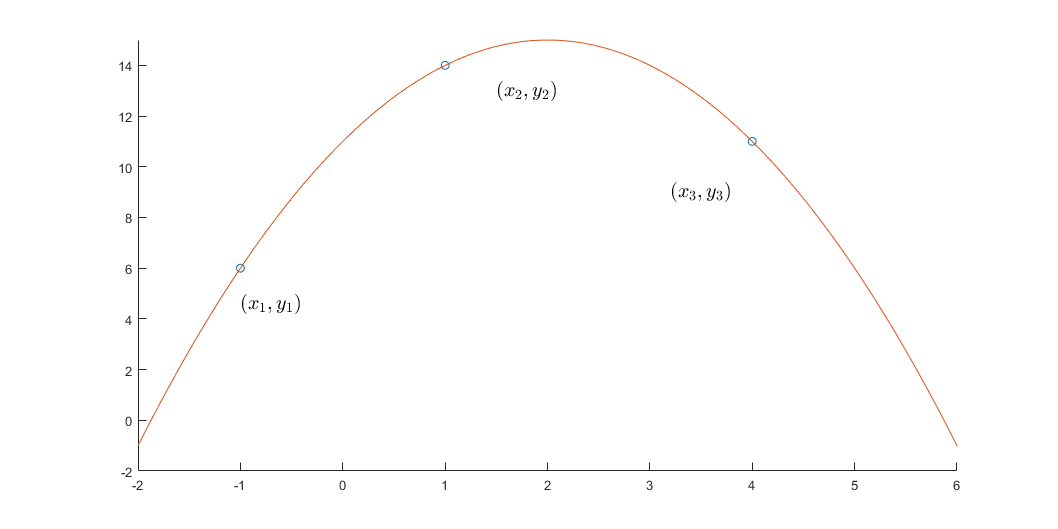}
		\caption{The polynomial interpolation of three points}
	\end{figure}
	
	We would construct three curves, where
	\begin{itemize}
		\item[$\bullet$] The first curve $f_1(x)$ takes value $1$ at $x_1$, while taking $0$ at other two points;
		\item[$\bullet$] The second curve $f_2(x)$ takes value $1$ at $x_2$, while taking $0$ at other two points;
		\item[$\bullet$] The third curve $f_3(x)$ takes value $1$ at $x_3$, while taking $0$ at other two points.
	\end{itemize}
	
	with the figures shown below.
	\clearpage
	
	\begin{figure}[h!] 
		\centering
		\includegraphics[width=14cm]{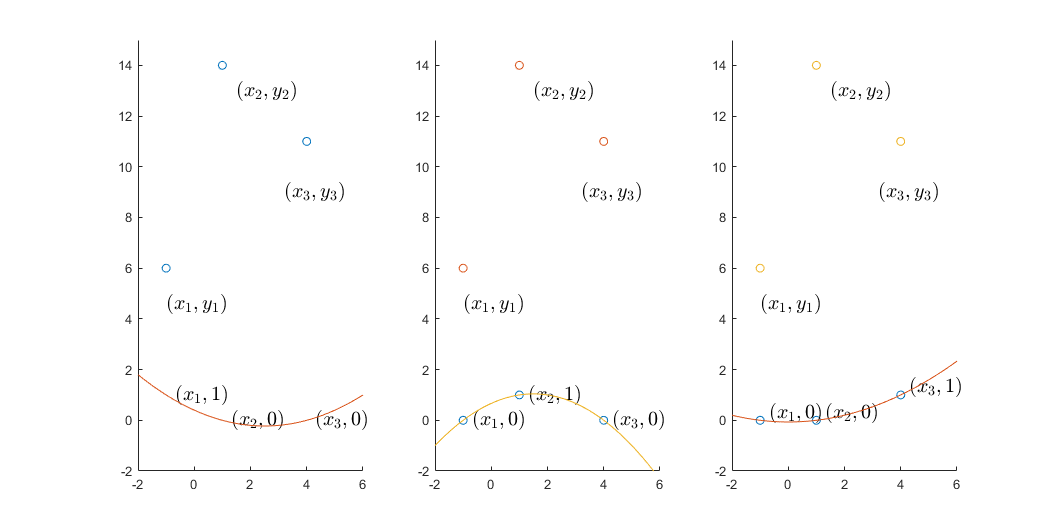}
		\caption{Three curves $f_i(x),i=1,2,3$}
	\end{figure}
	
	Hence
	
	\begin{itemize}
		\item[$\bullet$] $y_1f_1(x)$ can make sure take value $y_1$ at $x_1$, 0 at other points;
		\item[$\bullet$] $y_2f_2(x)$ can make sure take value $y_2$ at $x_2$, 0 at other points;
		\item[$\bullet$] $y_3f_3(x)$ can make sure take value $y_3$ at $x_3$, 0 at other points.
	\end{itemize}
	Hence, $f(x)=y_1f_1(x)+y_2f_2(x)+y_3f_3(x)$ must go through $(x_1,y_1),(x_2,y_2),(x_3,y_3)$.We can test and verify it by MATLAB:
	\begin{figure}[h!]
		\centering
		\includegraphics[width=14cm]{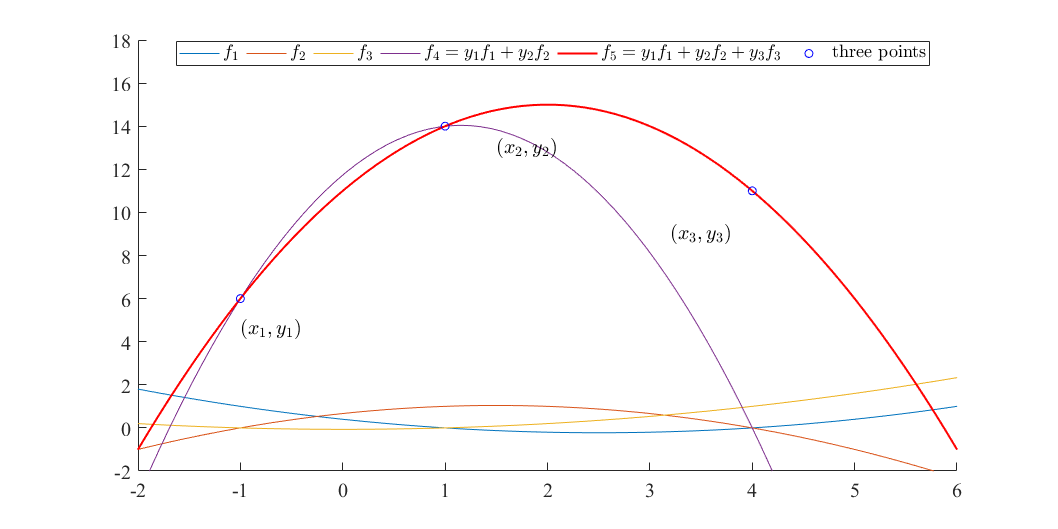} 
		\caption{The result of interpolation}
	\end{figure} 
	
	However, this is a simple case, so we need to standardize it in the next chapter so that it could be used more universally.
	
	\subsubsection{Standardization}
	The only work left is to standardize the procedure above. We use symbol $l_i(x_j)$,where $i=1,2,3\cdots n, j=1,2,3\cdots n$
	which satisfies:
	\begin{itemize}
		\item[$\bullet$] $l_i(x)$ must be a $n$-th order function, where $n$ is the number of data points minus $1$, considering we start from $0$.
		\item[$\bullet$] 
		$l_i(x)=
		\begin{cases}
			0& i=j,\\
			1& i\neq j.
		\end{cases}$
	\end{itemize}
	Hence we can construct $l_i(x_j)$:
	$$
	l_{i}(x)=\prod_{j \neq i}^{1 \leq j \leq n} \frac{\left(x-x_{j}\right)}{\left(x_{i}-x_{j}\right)}.
	$$
	Obviously, $l_i(x_j)$ satisfies the above conditions.
	
	So finally, we get
	$$
	L(x)=\sum_{i=1}^{n} y_{i} l_{i}(x),
	$$
	which is the \emph{Lagrange interpolation}, which could also be defined as follows:

	\begin{defn}
		The $n$-th order Lagrange interpolation for function $f(x)$ is
		$$
		L_{n}(x)=\sum_{i=0}^{n} \frac{\left(x-x_{1}\right)\left(x-x_{2}\right) \cdots\left(x-x_{i-1}\right)\left(x-x_{i+1}\right) \cdots\left(x-x_{n}\right)}{\left(x_{0}-x_{1}\right)\left(x_{0}-x_{2}\right) \cdots\left(x_{0}-x_{i-1}\right)\left(x_{0}-x_{i+1}\right) \cdots\left(x_{0}-x_{n}\right)} f\left(x_{i}\right),
		$$
		
		$i=0,1, \cdots, n$. Moreover, let
		$$
		\omega_{n+1}(x)=\prod_{j=0}^{n}\left(x-x_{j}\right),
		$$
		
		then
		$$
		\omega_{n+1}^{\prime}\left(x_{i}\right)=\prod_{j=0 \atop j \neq i}^{n}\left(x_{i}-x_{j}\right),
		$$
		
		then $L_{n}(x)$ can also be expressed as
		$$
		L_{n}(x)=\sum_{i=0}^{n} \frac{\omega_{n+1}(x)}{\left(x-x_{i}\right) \omega_{n+1}^{\prime}\left(x_{i}\right)} f\left(x_{i}\right)
		$$
	\end{defn}
	
	The \emph{Lagrange interpolation} is simple and trivial, but when we add one point, the whole process of calculation needs to be repeated. So we would introduce the \emph{Newton interpolation}.
	\subsection{Newton Interpolation}
	\subsubsection{Background}
	Constructing interpolation by the undetermined coefficient method and \emph{Lagrange interpolation} could solve problems in proper ways, but the workload of a computer would be unimaginable when it comes to a large amount of a number set. So the \emph{Newton interpolation} is introduced to overcome this shortage.
	\subsubsection{Main Idea}
	Suppose we have $n+1$ distinct points $(x_0,y_0),(x_1,y_1),\cdots,(x_n,y_n)$ and try to find function $f$ to interpolate these points such that $y_i=f(x_i), i=0,2,\cdots,n$.
	
	First we consider $f_1(x)$ passing through points $\big(x_0,f(x_0)\big), \big(x_1,f(x_1)\big)$:
	$$
	f_1(x)=f(x_0)+k_1(x-x_0).
	$$
	
	Note that the structure of $k(x-x_0)$ ensures $f_1(x_0)=f(x_0)$. By substituting $f_1(x_1)=f(x_1)$, we get the coefficient $k$.\\
	$$
	k_1=\frac{f\left(x_{1}\right)-f\left(x_{0}\right)}{x_{1}-x_{0}},
	$$
	
	hence
	$$
	f_1(x)=f(x_{0})+\frac{f(x_1)-f(x_{0})}{x_1-x_0}(x-x_0).
	$$
	
	With the same procedure, we can find the function $f_2$ for $\big((x_0,f_1(x_0)),(x_1,f(x_1)),(x_2,f(x_2)\big)$:
	$$
	f_2(x)=f_1(x)+k_2(x-x_0)(x-x_1).
	$$
	
	By solving $f_2(x_2)=f(x_2)$, we get:
	$$
	k_2=\frac{\frac{f\left(x_{2}\right)-f\left(x_{1}\right)}{x_{2}-x_{1}}-\frac{f\left(x_{1}\right)-f\left(x_{0}\right)}{x_{1}-x_{0}}}{x_{2}-x_{0}},
	$$
	hence
	$$
	\begin{aligned}
		&f_2(x)=f\left(x_{0}\right)+\frac{f\left(x_{1}\right)-f\left(x_{0}\right)}{x_{1}-x_{0}}\left(x-x_{0}\right)+\frac{\frac{f\left(x_{2}\right)-f\left(x_{1}\right)}{x_{2}-x_{1}}-\frac{f\left(x_{1}\right)-f\left(x_{0}\right)}{x_{1}-x_{0}}}{x_{2}-x_{0}}\left(x-x_{0}\right)\left(x-x_{1}\right).
	\end{aligned}
	$$
	Analysing the forms of $k_1$ and $k_2$, some regularities can be found, so we would introduce the \textit{difference quotient}.
	
	\begin{defn}
		The $k$-th order difference quotient at $x_{0}, x_{1}, \cdots, x_{k}$ of $f(x)$ is defined as follows:
		\begin{itemize}
			\item[$\bullet$] $f[x_i,x_j]=\dfrac{f(x_i)-f(x_j)}{x_i-x_j},i \neq j$ is the first order difference quotient;\\
			\item [$\bullet$] $f[x_i,x_j,x_k]=\dfrac{f[i,j]-f[j,k]}{x_i-x_k},i \neq j \neq k$ is the second order difference quotient.
		\end{itemize}
		And universally
		$$
		f\left[x_{0}, x_{1}, \cdots, x_{k}\right]=\frac{f\left[x_{1}, \cdots, x_{k}\right]-f\left[x_{0}, \cdots, x_{k-1}\right]}{x_{k}-x_{0}}
		$$
		is the $k$-th order difference quotient of $f(x)$.
	\end{defn}
	
	Hence we could use the \emph{difference quotients} to define and simplify the \emph{Newton interpolation}.
	\begin{defn}
		The Newton interpolation of points $\big(x_0,f(x_0)),(x_1,f(x_1)),...,(x_n,f(x_n)\big)$ is
		$$
		\begin{aligned}
			f(x)=& f\left(x_{0}\right)+f\left[x_{0}, x_{1}\right]\left(x-x_{0}\right)+f\left[x_{0}, x_{1}, x_{2}\right]\left(x-x_{0}\right)\left(x-x_{1}\right)+\cdots \\
			&+f\left[x_{0}, x_{1}, \cdots, x_{n-2}, x_{n-1}\right]\left(x-x_{0}\right)\left(x-x_{1}\right) \cdots\left(x-x_{n-2}\right)\left(x-x_{n-1}\right) \\
			&+f\left[x_{0}, x_{1}, \cdots, x_{n-1}, x_{n}\right]\left(x-x_{0}\right)\left(x-x_{1}\right) \cdots\left(x-x_{n-1}\right)\left(x-x_{n}\right).
		\end{aligned}
		$$
	\end{defn}
	
	We can follow the table to calculate the \emph{difference quotient} in \emph{Newton interpolation}, we can also see that by adding a new point, we only need to calculate only one new item instead of recalculating the whole procedure again.
	\begin{figure}[h!] 
		\centering
		\includegraphics[width=10.5cm]{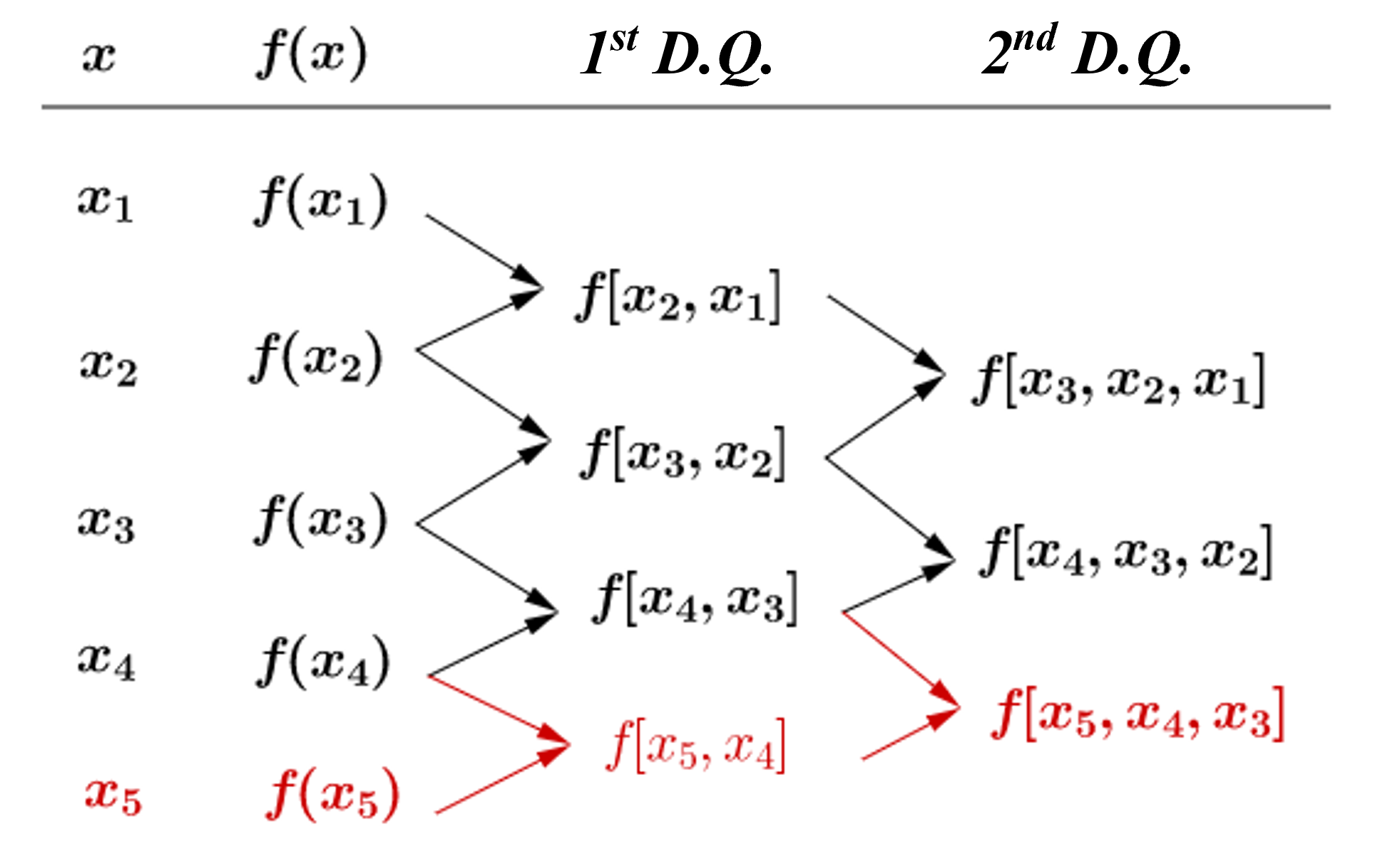}
		\caption{Table for \emph{difference quotients}}
	\end{figure} 
	
	Now that we have introduced two kinds of interpolation, it is necessary to point out that these two formulas form the same polynomial interpolation given the same points, which would be proved in the next section.
	
	\subsection{Uniqueness and Error}
	\subsubsection{Uniqueness}
	We have shown multiple methods to calculate the polynomial interpolation given numbers of points. The uniqueness theorem shows that the result is always the same if using different methods, which could be proved by linear algebra.
	
	Consider the case of constructing the polynomial interpolation by the undetermined coefficient method:
	$$
	\left\{\begin{array}{c}
		a_{0}+a_{1} x_{0}+\cdots a_{n} x_{0}^{n}=f\left(x_{0}\right) \\
		a_{0}+a_{1} x_{1}+\cdots a_{n} x_{1}^{n}=f\left(x_{1}\right) \\
		\cdots \\
		a_{0}+a_{1} x_{n}+\cdots a_{n} x_{n}^{n}=f\left(x_{n}\right).
	\end{array}\right.
	$$
	
	Usually, we solve a system of equations by converting it to the matrix form, so the \emph{Vandermonde matrix} would be introduced.
	
	\begin{defn} The $(n+1)$-th order Vandermonde matrix is
		$$
		\boldsymbol{V}=\left|\begin{array}{ccccc}
			1 & x_{0} & x_{0}^{2} & \cdots & x_{0}^{n} \\
			1 & x_{1} & x_{1}^{2} & \cdots & x_{1}^{n} \\
			\vdots & \vdots & \vdots & & \vdots \\
			1 & x_{n} & x_{n}^{2} & \cdots & x_{n}^{n}.
		\end{array}\right|
		$$
	\end{defn}
	
	Then the polynomial could be written in the matrix form, i.e.,
	$$
	\left(\begin{array}{l}
		f(x_0) \\
		f(x_1) \\
		\cdots \\
		f(x_n)
	\end{array}\right)
	=
	\boldsymbol{V}
	\left(\begin{array}{l}
		a_0 \\
		a_1 \\
		\cdots \\
		a_n
	\end{array}\right).
	$$
	
	We have learned from linear algebra, that if the whole system of equations has a unique solution if and only if the determinant of the \emph{Vandermonde matrix} is not $0$. When $i \neq j, x_{i} \neq x_{j}$, we have
	$$
	\left|\boldsymbol{V}\right|=\prod_{0 \leq i<j \leq n}\left(x_{j}-x_{i}\right) \neq 0,
	$$
	which means that when $x_{i}, i=0,1,\cdots,n$ are distinct points, the interpolation polynomial exists and is unique.
	
	Hence we can conclude that results of the undetermined coefficient method, \emph{Lagrange interpolation} and \emph{Newton interpolation} are the same.
	\subsubsection{Error}
	\begin{defn}
		Let $f(x)$ be the function passing through $n+1$ points $(x_0,y_0),(x_1,y_1),\cdots,$ $(x_n,y_n)$, and $P_n(x)$ be the $n$-th order interpolation, then the $n$-th order error is expressed as $R_n(x)$, and
		$$
		R_{n}(x)=f(x)-P_{n}(x)=\frac{f^{(n+1)}(\xi)}{(n+1) !} \omega_{n+1}(x),
		$$
		where
		$$
		\omega_{n+1}(x)=\prod_{j=0}^{n}\left(x-x_{j}\right).
		$$
	\end{defn}
	
	\textbf{Note:} This error could be applied for both \emph{Lagrange} and \emph{Newton interpolation} since they are the same expression, by the proof of the uniqueness theorem above.
	\clearpage
	\section{Neville's Algorithm——Stepwise Linear Interpolation}
	\subsection{Main Idea}
	If we know two points, the interpolation is to find a point on the straight line connected by two points naturally.( \textbf{Note}: That's the reason why we call it '\textbf{linear}' ) It can be verified that the following formula is the line from $(x_0, y_0)$ to $(x_1, y_1)\ (x_0\neq x_1)$:
	$$
	P_{0,1}(x)=\frac{x-x_{1}}{x_{0}-x_{1}} y_{0}+\frac{x-x_{0}}{x_{1}-x_{0}} y_{1},
	$$
	Then we consider the case of three points. We try to interpolate one point between $(x_0,y_0)$ and $(x_1,y_1)$, interpolate another point between $(x_1,y_1)$ and $(x_2,y_2)$:
	$$
	P_{1,2}(x)=\frac{x-x_{2}}{x_{1}-x_{2}} y_{1}+\frac{x-x_{1}}{x_{2}-x_{1}} y_{2}.
	$$
	We acquiesce that two points are very closed to $(x_0,y_0),(x_1,y_1)$
	$$
	\frac{x-x_{2}}{x_{0}-x_{2}} P_{0,1}+\frac{x-x_{0}}{x_{2}-x_{0}} P_{1,2}=P_{0,1,2}
	$$
	By elimination we can get:
	$$
	\frac{\left(x-x_{1}\right)\left(x-x_{2}\right)}{\left(x_{0}-x_{1}\right)\left(x_{0}-x_{2}\right)} y_{0}+\frac{\left(x-x_{0}\right)\left(x-x_{2}\right)}{\left(x_{1}-x_{0}\right)\left(x_{1}-x_{2}\right)} y_{1}+\frac{\left(x-x_{0}\right)\left(x-x_{1}\right)}{\left(x_{2}-x_{0}\right)\left(x_{2}-x_{1}\right)} y_{2}=P_{0,1,2},
	$$
	so we could find the recursion:
	$$
	P_{0,1 \cdots, n}=\frac{x-x_{0}}{x_{n}-x_{0}} P_{1,2 \cdots, n}(x)+\frac{x-x_{n}}{x_{0}-x_{n}} P_{0,1,2 \cdots, n-1}(x).
	$$
	
	We can generalize the steps as the above figure:
	\begin{figure}[htbp] 
		\centering
		\includegraphics[width=7.5cm]{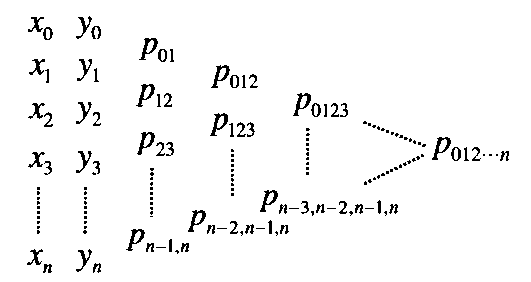}
		\caption{Recursion for the Neville's algorithm}
	\end{figure} 
	
	We could interpolate step by step, so that's the reason why we call it \textbf{'Stepwise'}.
	
	Hence we could get the theorem as follows.
	\begin{theorem}[\citet{cite-ex2}]
		Let a function $f$ be defined at points $x_{0}, x_{1}, \cdots, x_{k}$ where $x_{j}$ and $x_{i}$ are two distinct members. For each $k$, there exists a Lagrange polynomial $P$ that interpolates the function $f$ at the $k+1$ points $x_{0}, x_{1}, \cdots, x_{k} .$ The $k$ th Lagrange polynomial is defined as:
		$$
		P(x)=\frac{\left(x-x_{j}\right) P_{0,1, \cdots, j-1, j+1, \cdots, k}(x)-\left(x-x_{i}\right) P_{0,1, \cdots, i-1, i+1, \cdots, k}(x)}{\left(x_{i}-x_{j}\right)}.
		$$
		
		The $P_{0,1, \cdots, j-1, j+1, \cdots, k}$ and $P_{0,1, \cdots, i-1, i+1, \cdots, k}$ are often denoted $\hat{Q}$ and $Q$, respectively, for ease of notation.
		$$
		P(x)=\frac{\left(x-x_{j}\right) \hat{Q}(x)-\left(x-x_{i}\right) Q(x)}{\left(x_{i}-x_{j}\right)}.
		$$
	\end{theorem}
	
	The procedure that uses the result of Theorem 2.1 to recursively generate interpolating polynomial approximations is called \emph{Neville’s method}.\citep{cite-ex2}
	\subsection{Advantages and Disadvantages}
	\subsubsection{Advantages}
	The \textbf{advantage} of \emph{Neville's interpolation} is that if the data are arranged in order 
	of closeness to the interpolated point, is that none of the work performed to obtain a specific degree result 
	must be redone to evaluate the next higher degree result.\citep{cite-ex3}
	\subsubsection{Disadvantages}
	The \textbf{disadvantage} of \emph{Neville's interpolation} is that all of the work must be redone for each new value of $x$. The amount of work is essentially the same as for a Lagrange polynomial. The divided difference polynomial minimizes these disadvantages.\citep{cite-ex3}
	\clearpage
	\section{Hermite Interpolation: About Derivatives}
	\subsection{Main Idea}
	We have introduced some methods of interpolation, which usually concerns the value of points, say, $y_i$. However, in practice, the derivatives of points, say, $y'_i$, are also important since they represent the rates of change and we do not wish to discard this information. That is the reason we introduce the \emph{Hermite interpolation}, which  requires only one more rule than the \emph{Lagrange interpolation}: The $n$-th order interpolation function $H_i(x)$ has the same derivative with $y_i$, i.e.,
	
	$$
	H'_i(x)=y'_i, i=0,1,\cdots,n.
	$$
	
	where $n$ is the number of conditions minus $1$, considering we start from $0$.
	
	We can compare the \emph{Hermite interpolation} with the main idea of the \emph{Lagrange interpolation} —— to construct 'base functions'. Without loss of generality, we discuss a two points, $3$-order case of the \emph{Hermite interpolation}. Usually, when we want to have a $n$-th order interpolation, we may need $n$ distinct points, but for the \emph{Hermite interpolation}, we only need almost half of the points, since we have twice conditions than before: the additional derivative conditions. Suppose we have two points $(x_0,y_0),(x_1,y_1)$, when we construct a function $H_3(x)$, it needs to satisfy the following conditions: 
	$$
	\begin{cases}
		H_3(x_i)&=y_i\\ 
		H'_3(x_i)&=y'_i
	\end{cases}
	$$
	
	where $i=0,1$.
	
	Just using the undetermined coefficient method, let $H_3(x)=ax^3+bx^2+cx+d$. It would be hard to calculate and not easy to generalize into higher orders, so we introduce the method of 'base function': note four base functions $\alpha_0(x),\alpha_1(x),\beta_0(x),\beta_1(x)$, which satisfy:
	$$
	\left\{\begin{array} { l } 
		{ \alpha _ { 0 } ( x _ { 0 } ) = 1 }; \\
		{ \alpha _ { 0 } ( x _ { 1 } ) = 0 }; \\
		{ \alpha _ { 0 } ^ { \prime } ( x _ { 0 } ) = 0 }; \\
		{ \alpha _ { 0 } ^ { \prime } ( x _ { 1 } ) = 0 }.
	\end{array} \quad \left\{\begin{array} { l } 
		{ \alpha _ { 1 } ( x _ { 0 } ) = 0 }; \\
		{ \alpha _ { 1 } ( x _ { 1 } ) = 1 }; \\
		{ \alpha _ { 1 } ^ { \prime } ( x _ { 0 } ) = 0 }; \\
		{ \alpha _ { 1 } ^ { \prime } ( x _ { 1 } ) = 0 }.
	\end{array} \quad \left\{\begin{array} { l } 
		{ \beta _ { 0 } ( x _ { 0 } ) = 0 }; \\
		{ \beta _ { 0 } ( x _ { 1 } ) = 0 }; \\
		{ \beta _ { 0 } ^ { \prime } ( x _ { 0 } ) = 1 }; \\
		{ \beta _ { 0 } ^ { \prime } ( x _ { 1 } ) = 0 }.
	\end{array} \quad \left\{\begin{array}{l}
		\beta_{1}\left(x_{0}\right)=0; \\
		\beta_{1}\left(x_{1}\right)=0; \\
		\beta_{1}^{\prime}\left(x_{0}\right)=0; \\
		\beta_{1}^{\prime}\left(x_{1}\right)=1.
	\end{array}\right.\right.\right.\right. .
	$$ 
	Note that $\alpha_0(x),\alpha_1(x),\beta_0(x),\beta_1(x)$ are all $3$-order functions.
	Let $H_3=y_0 \alpha_0(x)+y_1\alpha_1(x)+y'_0\alpha_0(x)+y_1\alpha_1(x)$, such that $H_3$ is a polynomial interpolation which has an order less than three.
	
	The remaining problem is how to find the base function. First, we try to find $\alpha_0(x)$, from the given condition we have:
	\begin{equation*}
		\begin{cases}
			\begin{cases}
				\alpha_{0}\left(x_{1}\right)&=0 \\
				\alpha_{0}^{\prime}\left(x_{1}\right)&=0
			\end{cases}
			&\Rightarrow
			\alpha_{0}(x)=\left[a+b\left(x-x_{0}\right)\right]\left(x-x_{1}\right)^{2},\\
			\alpha_{0}\left(x_{0}\right)=1 &\Rightarrow a=\dfrac{1}{\left(x_{0}-x_{1})^2\right.}\\
			{\alpha}_{0}^{\prime}\left({x}_{0}\right)=0, &\Rightarrow b=\dfrac{2}{\left(x_{1}-x_{0}\right)\left(x_{0}-x_{1}\right)^{2}}.
		\end{cases}
	\end{equation*}
	From above all, we can get:
	$$
	\alpha_{0}(x)=\left(1+2 \frac{x-x_{0}}{x_{1}-x_{0}}\right)\left(\frac{x-x_{1}}{x_{0}-x_{1}}\right)^{2}.
	$$
	Similarly, we get
	\begin{align*}
		\begin{cases}
			\alpha_{1}(x)&=\left[1+\dfrac{2\left(x-x_{0}\right)}{x_{1}-x_{0}}\right]\left(\dfrac{x-x_{0}}{x_{1}-x_{0}}\right)^{2},\\
			\beta_{0}(x)&=\left(x-x_{0}\right)\left(\dfrac{x-x_{1}}{x_{0}-x_{1}}\right)^{2},\\
			\beta_{1}(x)&=\left(x-x_{1}\right)\left(\dfrac{x-x_{0}}{x_{1}-x_{0}}\right)^{2}.\\
		\end{cases}
	\end{align*}
	So the $3$-order $Hermite \ Interpolation$ can be expressed as:
	$$
	\begin{aligned}
		H_{3}(x)=&\left[\left(1+2 \frac{x-x_{0}}{x_{1}-x_{0}}\right) y_{0}+\left(x-x_{0}\right) y_{0}^{\prime}\right]\left(\frac{x-x_{1}}{x_{0}-x_{1}}\right)^{2}\\
		&+\left[\left(1+2 \frac{x-x_{1}}{x_{0}-x_{1}}\right) y_{1}+\left(x-x_{1}\right) y_{1}^{\prime}\right]\left(\frac{x-x_{0}}{x_{1}-x_{0}}\right)^{2}.
	\end{aligned}
	$$
	\subsection{Standardization}
	\begin{defn}[\citet{cite-ex2}]
		If $f \in C^{1}[a, b]$ and $x_{0}, \ldots, x_{n} \in[a, b]$ are distinct, the unique polynomial of least degree agreeing with $f$ and $f^{\prime}$ at $x_{0}, \ldots, x_{n}$ is the Hermite polynomial of degree at most $2 n+1$ given by
		$$
		H_{2 n+1}(x)=\sum_{j=0}^{n} f\left(x_{j}\right) H_{n, j}(x)+\sum_{j=0}^{n} f^{\prime}\left(x_{j}\right) \hat{H}_{n, j}(x)
		$$
		
		where, for $L_{n, j}(x)$ denote the $j$-th Lagrange coefficient polynomial of degree $n$, we have
		$$
		H_{n, j}(x)=\left[1-2\left(x-x_{j}\right) L_{n, j}^{\prime}\left(x_{j}\right)\right] L_{n, j}^{2}(x) \quad \text { and } \quad \hat{H}_{n, j}(x)=\left(x-x_{j}\right) L_{n, j}^{2}(x) .
		$$
	\end{defn}
	
	\subsection{Error}
	\begin{defn}[\citet{cite-ex4}]
		If $H_{2n+1}(x)$ is the Hermite interpolation polynomial of $f(x)$ at points   $x_0,x_1,\dots,x_n (x_0,x_1,\dots,x_n \in [a,b])$, then the $(2n+1)$-th order error is expressed as $R_{2n+1}(x)$, and
		$$
		R_{2n+1}(x)=f(x)-H_{2 n+1}(x)=\frac{f^{(2 n+2)}(\xi_x)}{(2 n+2) !} \omega^{2}(x),
		$$
		
		for some $\xi_x \in(a, b)$, and
		$$
		\omega(x)=\prod_{i=0}^n(x-x_i).
		$$
	\end{defn}
	
	\clearpage
	\section{About Runge's Phenomenon - Cubic Splines and Chebyshev Nodes}
	We have shown above that if we have $n+1$ points $(x_0,y_0), (x_2,y_2), \cdots, (x_n,y_n)$, then we can have a $n$-th order polynomial by any kind of interpolation. However, there would be some problems encountered when $n$ is large.
	\subsection{Runge's Phenomenon}
	It could be summarized that \emph{Runge phenomenon} is an oscillation problem that occurs at the edges of an interval containing data points when using polynomial interpolation with a high degree of a polynomial over a set of equispaced interpolation points. 
	
	\emph{Runge's phenomenon} states that when the order of the polynomial is too high, it would cause high errors, which might diverge to infinity.
	
	\begin{figure}[h!]
		\centering
		\includegraphics[width=14cm]{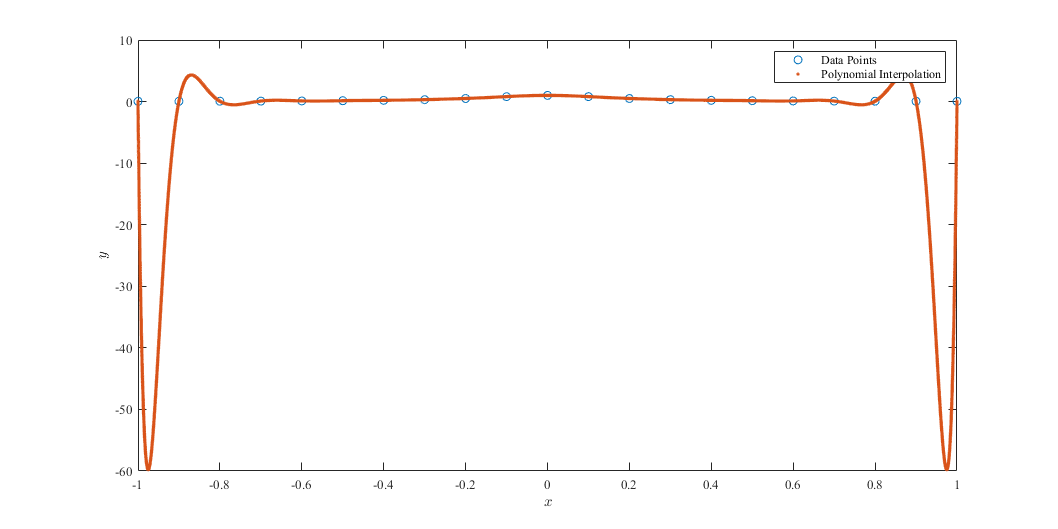}
		\caption{Demonstration of Runge's phenomenon for $y=\frac{1}{1+25x^2}, x\in[-1,1]$}
	\end{figure}
	
	\subsubsection{Proof}
	The Weierstrass approximation theorem states that for every continuous function $f(x)$ defined on an interval $[a,b]$, there exists a set of polynomial functions $P_n(x)$ for $n=0, 1, 2,\cdots$, each of degree at most $n$, that approximates $f(x)$ with uniform convergence over $[a,b]$ as $n$ tends to infinity, that is,
	
	$$
	\lim _{n \rightarrow \infty}\left(\max _{a \leq x \leq b}\left|f(x)-P_{n}(x)\right|\right)=0.
	$$
	
	Consider the case where one desires to interpolate through $n+1$ equispaced points of a function $f(x)$ using the n-degree polynomial $P_n(x)$ that passes through those points. Naturally, one might expect from Weierstrass' theorem that using more points would lead to a more accurate reconstruction of $f(x)$. However, this particular set of polynomial functions $P_n(x)$ is not guaranteed to have the property of uniform convergence; the theorem only states that a set of polynomial functions exists, without providing a general method of finding one.
	
	The $P_n(x)$ produced in this manner may diverge away from $f(x)$ as n increases; this typically occurs in an oscillating pattern that magnifies near the ends of the interpolation points. This phenomenon is attributed to Runge.\citep{cite-ex1}
	
	\subsection{Solution 1: Piecewise Polynomials}
	Instead of using the $n$-th order interpolation polynomial considering $n$ is too large, we would use multiple polynomials among small intervals, one common application is the \emph{Cubic spline}, which is to use multi $3$-order polynomials and combine them to the final interpolation formula instead.
	
	\begin{defn}{\citep{cite-ex2}}
		Given a function $f$ defined on $[a, b]$ and a set of nodes $a=x_{0}<x_{1}<\cdots<$ $x_{n}=b$, a cubic spline interpolant $S$ for $f$ is a function that satisfies the following conditions:
		\begin{description}
			\item[(a)] $S(x)$ is a cubic polynomial, denoted $S_{j}(x)$, on the subinterval $\left[x_{j}, x_{j+1}\right]$ for each $j=0,1, \ldots, n-1$;
			\item[(b)] $S_{j}\left(x_{j}\right)=f\left(x_{j}\right)$ and $S_{j}\left(x_{j+1}\right)=f\left(x_{j+1}\right)$ for each $j=0,1, \ldots, n-1$;
			\item[(c)] $S_{j+1}\left(x_{j+1}\right)=S_{j}\left(x_{j+1}\right)$ for each $j=0,1, \ldots, n-2$ (Implied by \textbf{(b)});
			\item[(d)] $S_{j+1}^{\prime}\left(x_{j+1}\right)=S_{j}^{\prime}\left(x_{j+1}\right)$ for each $j=0,1, \ldots, n-2$;
			\item[(e)] $S_{j+1}^{\prime \prime}\left(x_{j+1}\right)=S_{j}^{\prime \prime}\left(x_{j+1}\right)$ for each $j=0,1, \ldots, n-2$;
			\item[(f)] One of the following sets of boundary conditions is satisfied:
			\begin{description}
				\item[(i)] $S^{\prime \prime}\left(x_{0}\right)=S^{\prime \prime}\left(x_{n}\right)=0 \quad$ (\textbf{natural} (or \textbf{free}) boundary);
				\item[(ii)] $S^{\prime}\left(x_{0}\right)=f^{\prime}\left(x_{0}\right)$ and $S^{\prime}\left(x_{n}\right)=f^{\prime}\left(x_{n}\right) \quad$ (\textbf{clamped} boundary).
			\end{description}
		\end{description}
	\end{defn}
	
	When the free boundary conditions occur, the spline is called a \textbf{natural spline}\citep{cite-ex2}. Note that the \textbf{natural spline} is the most commonly used spline.
	
	\begin{theorem}{\citep{cite-ex2}}
		If $f$ is defined at $a=x_{0}<x_{1}<\cdots<x_{n}=b$, then $f$ has a unique natural spline interpolant $S$ on the nodes $x_{0}, x_{1}, \ldots, x_{n}$; that is, a spline interpolant that satisfies the natural boundary conditions $S^{\prime \prime}(a)=0$ and $S^{\prime \prime}(b)=0$.
	\end{theorem}
	
	\subsubsection{Construction}
	We could always get a system of $4n$ equations if we are given enough information which is stated in the definition above, and we always try to solve by algorithms and apply them in computers when $n$ is large. Here follows the algorithm for calculating the natural cubic spline, note that the time complexity is $\mathcal{O}(n)$.
	
	%\clearpage
	\begin{algorithm}[h!]
		\caption{Algorithm for natural cubic spline\citep{cite-ex2}} % 算法的名字
		\hspace*{0.02in} {\bf Input:} % 算法的输入， \hspace*{0.02in} 用来控制位置，同时利用 \\ 进行换行
		$n ; x_{0}, x_{1}, \ldots, x_{n} ; a_{0}=f\left(x_{0}\right), a_{1}=f\left(x_{1}\right), \ldots, a_{n}=f\left(x_{n}\right)$\\
		\hspace*{0.02in} {\bf Output:} % 算法的结果输出
		$a_{j}, b_{j}, c_{j}, d_{j} \text { for } j=0,1, \ldots, n-1$
		\begin{algorithmic}[1]
			\For{$i=0,1, \ldots, n-1$} % For 语句，需要和 EndFor 对应
			\State set $h_{i}=x_{i+1}-x_{i}$.
			\EndFor
			\For{$i=1,2,\cdots,n-1$}
			\State set $\alpha_{i}=\frac{3}{h_{i}}\left(a_{i+1}-a_{i}\right)-\frac{3}{h_{i-1}}\left(a_{i}-a_{i-1}\right)$.
			\EndFor
			\State Set $l_0=1$;\\
			$\quad \ \ \mu_0=0$;\\
			
			$\quad \ \ z_0=0$.
			\For{$i=1,2,\cdots,n-1$}
			\State set 
			$l_{i}=2\left(x_{i+1}-x_{i-1}\right)-h_{i-1} \mu_{i-1}$;\\
			
			$\quad \quad \quad \mu_{i}=h_{i} / l_{i}$;\\
			
			$\quad \quad \quad z_{i}=\left(\alpha_{i}-h_{i-1} z_{i-1}\right) / l_{i}$.
			\EndFor
			\State Set
			$l_{n}=1$;\\
			
			$\quad \ \ z_{n}=0$;\\
			
			$\quad \ \ c_{n}=0$.
			
			\For{$j=n-1,n-2,\cdots,0$}
			\State set $c_{j}=z_{j}-\mu_{j} c_{j+1}$;\\
			
			$\quad \quad \quad b_{j}=\left(a_{j+1}-a_{j}\right) / h_{j}-h_{j}\left(c_{j+1}+2 c_{j}\right) / 3$;\\
			
			$\quad \quad \quad d_{j}=\left(c_{j+1}-c_{j}\right) /\left(3 h_{j}\right)$.
			\EndFor
			
			\State \Return $a_j,b_j,c_j,d_j$ for $0,1,\cdots,n-1$
		\end{algorithmic}
	\end{algorithm}
	
	\subsection{Solution 2: Chebyshev Nodes}
	We have stated that the \emph{Runge's phenomenon} would happen if we choose a set of equispaced interpolation points, which means it could be avoided by selecting non-equispaced points, say, \emph{Chebyshev nodes}.
	
	\begin{defn}
		Let $z=e^{i \theta}$ be a point on the unit circle. The associated $x$ coordinate is $x=\cos \theta$ or $\theta=\cos ^{-1} x$ where $x \in[-1,1]$. Define the $n$th degree Chebyshev polynomial to be $T_{n}(x)=\cos n \theta$. The Chebyshev nodes $x_0,x_1,\cdots,x_n$ are the roots of Chebyshev polynomial $T_{n+1}$.
	\end{defn}
	
	The reason why Chebyshev nodes work should be originated from the error between the generating function and the interpolating polynomial of order $n$, which is given by
	$$
	f(x)-P_{n}(x)=\frac{f^{(n+1)}(\xi)}{(n+1) !} \prod_{i=1}^{n+1}\left(x-x_{i}\right)
	$$
	
	for some $\xi$ in $(-1,1)$. Thus,
	$$
	\max _{-1 \leq x \leq 1}\left|f(x)-P_{n}(x)\right| \leq \max _{-1 \leq x \leq 1} \frac{\left|f^{(n+1)}(x)\right|}{(n+1) !} \max _{-1 \leq x \leq 1} \prod_{i=0}^{n}\left|x-x_{i}\right|.
	$$
	
	For the case of the Runge function which is interpolated at equidistant points, each of the two multipliers in the upper bound for the approximation error grows to infinity with $n$. Since $\frac{f^{(n+1)}(x)}{(n+1)!}$ could be regarded as a constant, we need to optimize $\prod_{i=0}^{n}\left|x-x_{i}\right|$ instead. If we choose \emph{Chebyshev nodes} from $x_0$ to $x_n$, then
	$$
	\left(x-x_{0}\right)\left(x-x_{1}\right) \ldots\left(x-x_{n}\right)=\frac{T_{n+1}(x)}{2^{n}}
	$$
	
	where
	$$
	x_k=\cos \left[\frac{(2 k+1) \pi}{2(n+1)}\right] \quad k=0,1, \ldots, \mathrm{n}.
	$$
	
	Note that $\frac{T_{n+1}(x)}{2^{n}}$ could have the smallest $\|\cdot\|_{\infty}$ value over the interval $[-1,1]$. We could prove it by assuming there exists $q_{n+1}$ which $\|\cdot\|_{\infty}$ value is not greater than $T_{n+1}$.
	
	Now $\left\|\frac{T_{n+1}}{2^{n}}\right\|_{\infty}=2^{-n}$ is achieved $n+2$ times within $[-1,1]$. By definition $|q_{n+1}(x)|<\dfrac{1}{2^n}$ at each $n+2$ extreme points. Thus 
	$$
	D(x)=\frac{T_{n+1}}{2^{n}}-q_{n+1}
	$$
	is a polynomial of degree $\leq n$ and has the same sign as $T_{n+1}$ at each of the $n+2$ extreme points.
	
	Note that $D(x)$ must change sign $n+1$ times on $[-1,1]$ which is impossible for a polynomial of degree $\leq n$, which would yield a contradiction. So we would get better results if using \emph{Chebyshev nodes} rather than \emph{equispaced nodes}.\citep{cite-ex8}
	
	Note that \emph{Chebyshev nodes} are only between $[-1,1]$, but we could rescale any intervals to $[-1,1]$ to achieve universality\citep{cite-ex2}.
	\begin{theorem}{\citep{cite-ex8}}
		For $t\in[-1,1]$, let
		$$
		x=\dfrac{t(b-a)+a+b}{2},
		$$
		
		then we can rescale any interval $[a,b]$ to $[-1,1]$.
	\end{theorem}
	\clearpage
	\section{Data Fitting}
	\subsection{Key Difference}
	In interpolation, we construct a curve through the data points. In doing so, we make the implicit assumption that the data points are accurate and distinct. Data fitting is applied to data that contain scatter, like noise, usually due to measurement errors. Here we want to find a smooth curve that approximates the data in some sense. Thus the curve does not necessarily hit the data points. \citep{cite-ex6}
	\begin{figure}[h!]
		\centering
		\includegraphics[width=12cm]{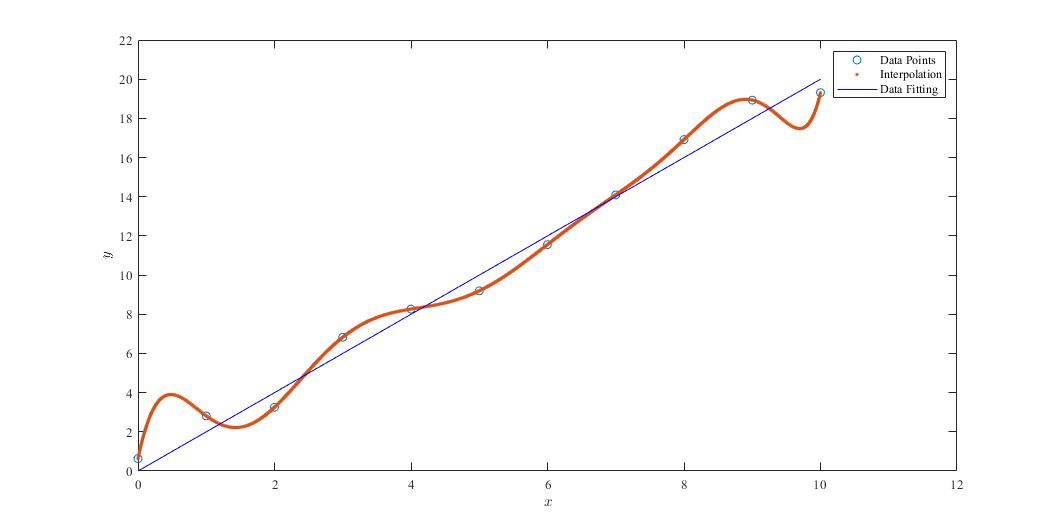}
		\caption{The figure of the comparison of interpolation and data fitting}
	\end{figure}

	\subsection{Linear Squares Method}
	\subsubsection{Main Idea}
	The \emph{least squares method} is introduced to decide the parameters of a curve that when the sum of distances squared from different points to the curve achieves the minimum.
	
	We consider an example here, suppose we have a series of data points as follows:
	\begin{table}[h!]
		\centering
		\begin{tabular}{|c|c|c|c|c|c|c|}
			\hline
			\textbf{$i$}   & \textbf{$1$}     & \textbf{$2$}     & \textbf{$3$}     & \textbf{$4$}    & \textbf{$5$}   & \textbf{$6$}    \\ \hline
			\textbf{$x_i$} & \textbf{$-9.19$} & \textbf{$-5.26$} & \textbf{$-1.39$} & \textbf{$6.71$} & \textbf{$4.7$} & \textbf{$2.66$} \\ \hline
			\textbf{$y_i$} & \textbf{$-8.01$} & \textbf{$6.78$}  & \textbf{$-1.47$} & \textbf{$4.71$} & \textbf{$4.1$} & \textbf{$4.23$} \\ \hline
		\end{tabular}
		\caption{A table of data points}
	\end{table}
	
	Now we want to find a line that satisfies the \emph{linear squares} condition.
	
	For our convenience, we consider only 3 points $A_1:(-9.19,-8.01),$  $A_2:(-5.26,6.78),$   $A_3:(-1.39,-1.47)$. Denote the distances from points $A_1, A_2, A_3$ to the line $y=ax+b$ as $d_1,d_2,d_3$ respectively. We want to minimize 
	$$
	D^2=d_1^2+d_2^2+d_3^2.
	$$
	
	We get
	\begin{equation*}
		\begin{cases}
			d_1&=-9.19a+b+8.01,\\
			d_2&=-5.26a+b-6.78,\\
			d_3&=-1.39a+b+1.47.
		\end{cases}
	\end{equation*}
	
	Hence:
	$$
	D^2=(-9.19a+b+8.01)^2+(-9.19a+b+8.01)^2+(-9.19a+b+8.01)^2.
	$$
	With the knowledge in calculus, we can get the solution of $a,b$ when $D^2$ achieves the minimum by calculating its derivative. The result for coefficients is
	
	\begin{equation*}
		\begin{cases}
			a&\approx 0.958,\\
			b&\approx -0.584.
		\end{cases}
	\end{equation*}
	
	So the fitting line by \emph{least squares method} of our example is :
	$$ 
	y=0.958x-0.584.
	$$
	
	\begin{figure}[h!]
		\centering
		\includegraphics[width=14cm]{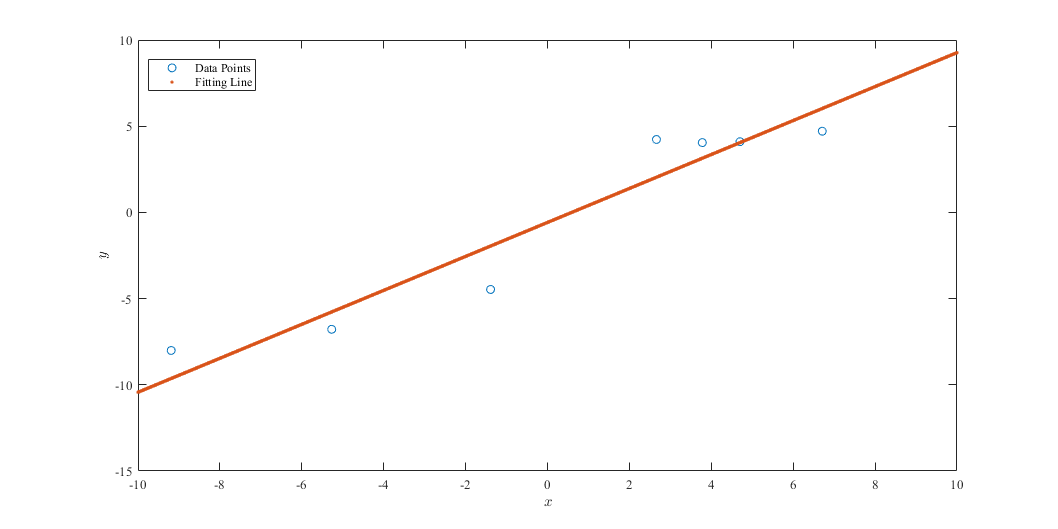}
		\caption{The figure of the fitting line}
	\end{figure} 
	
	\subsubsection{Advantages}
	The least squares method is the most convenient procedure for determining best linear approximations, but there are also important theoretical considerations that favour it. The minimax approach generally assigns too much weight to a bit of data that is badly in error, whereas the absolute deviation method does not give sufficient weight to a point that is considerably out of line with the approximation. The least squares approach puts substantially more weight on a point that is out of line with the rest of the data, but will not permit that point to completely dominate the approximation. An additional reason for considering the least squares approach involves the study of the statistical distribution of
	error.\citep{cite-ex7}
	
	\subsubsection{Some other forms}
	Sometimes we assume that the data exponentially distributed, which would require us to approximate the form
	$$
	y=be^{ax}
	$$
	
	or
	$$
	y=bx^a
	$$
	
	for $a,b\in \mathbb{C}$. So we need to minimize
	$$
	E=\sum_{i=1}^{n}\left(y_{i}-b e^{a x_{i}}\right)^{2},
	$$
	or
	$$
	E=\sum_{i=1}^{n}\left(y_{i}-b x_{i}^{a}\right)^{2}.
	$$
	
	\subsubsection{Standardization}
	Given function $f(x,\alpha_1,\alpha_2,\dots,\alpha_n)$, and $n$ distinct points $(x_1,y_1),(x_2,y_2),\dots,(x_n,y_n)$ on $f$. We would give the definition of \emph{residuals}:
	\begin{defn}
	The residual $r_i$ for function $f(x,\alpha_1,\alpha_2,\dots,\alpha_n)$ and $n$ distinct points $(x_1,y_1), $ $(x_2,y_2),\dots,(x_n,y_n)$ is
	$$
	r_i=f(x,\alpha_1,\alpha_2,\dots,\alpha_n)-y_i, \ i=1,2,...,n.
	$$
	\end{defn}

	The aim is to minimize the sum of residuals, that is :
	$$
	\min\sum_{i=1}^{N} r_{i}^{2}=\min\sum_{i=1}^{N}\left[f\left(x_{i} ; \alpha_{1}, \alpha_{2}, \cdots, \alpha_{n}\right)-y_{i}\right]^{2}
	$$
	
	which would simplify to \emph{normal equations}. This would be discussed in the later part.
	
	\subsection{Normal Equations}
	\subsubsection{Main Idea}
	Suppose we have $m$ samples, each sample has $n$ continuous features and a label value $y$. Now it is necessary to find out the linear relationship between features and label values. We can define such a \emph{cost function}:
	\begin{defn}
		The cost function for $m$ samples is
		$$
		\begin{aligned}
			\mathrm{J}(\mathrm{\theta_0,\theta_1,...,\theta_m}) &=\left(\mathrm{h_\theta}\left(x_{1}\right)-y_{1}\right)^{2}+\left(\mathrm{h_\theta}\left(x_{2}\right)-y_{2}\right)^{2}+\cdots+\left(\mathrm{h_\theta}\left(x_{m}\right)-y_{m}\right)^{2} \\
			&=\sum_{i=1}^{m}\left[\mathrm{~h_\theta}\left(x_{i}\right)-y_{i}\right]^{2}.
		\end{aligned}
		$$
	\end{defn}

	Note that
	\begin{itemize}
		\item [$\bullet$]$y_i$ is the true value of the $i$-th training sample;
		\item [$\bullet$]$h_\theta(x_1)$ is  prediction function of the $i$-th training sample;
		\item [$\bullet$] $\theta$ are the coefficients in the loss function.
	\end{itemize}
	
	Looking at the form of $J(\theta)$, we can relate it with $D^2$ we talked about before. We want to find coefficients set $\theta={\theta_1,\theta_2,\dots,\theta_n}$, which minimizes the value of loss function.We just need to find the partial derivative:
	$$
	\frac{\partial}{\partial \theta_{j}} J(\Theta)=\cdots \stackrel{\text { }}{=} 0, j=1,2,\cdots,n.
	$$
	
	By solving it we can get, by matrix form,
	$$
	\Theta=\left(\mathbf{X}^{T}\cdot\mathbf{X}\right)^{-1} \mathbf{X}^{T} \mathbf{Y}.
	$$
	We will get a detailed proof and calculation later. Now we look at an example of a loan record in a bank, where $x$ represents different features and $y$ represent the value:
	\begin{table}[h!]
		\centering
		\begin{tabular}{c|c|c|c|c|}
			$x_{0}$ & Income $x_{1}$ & Age $x_{2}$ & Credit level $x_{3}$ & Values $y$ \\
			\hline 1 & 3000 & 33 & 5 & 4500 \\
			1 & 6500 & 30 & 2 & 12000 \\
			1 & 3500 & 31 & 4 & 8500 \\
			1 & 2800 & 28 & 5 & 4000
		\end{tabular}
		\caption{Example: loan record in a bank}
	\end{table}
	
	Then we get
	$$
	\mathbf{X}=\left[\begin{array}{ccccc}
		1 & 3000 & 33& 5  \\
		1 & 6500 & 30 & 2  \\
		1 & 3500 & 31 & 4 \\
		1 & 2800 & 28 & 5 
	\end{array}\right] \quad \mathbf{y}=\left[\begin{array}{c}
		4500 \\
		12000 \\
		8500 \\
		4000
	\end{array}\right]
	$$
	By the method of \emph{normal equations}, we can calculate:
	
	$$
	\Theta = \left(\mathbf{X}^T\cdot \mathbf{X}\right)^{-1}\cdot \mathbf{X}^{T}\cdot \mathbf{y}.
	$$
	
	$$
	\begin{small}
	\Theta=
	\left(\left[\begin{array}{llll}
		1 & 1 & 1 & 1 \\
		3000 & 6500 & 3500 & 2800 \\
		33 & 30 & 31 & 28 \\
		5 & 2 & 4 & 5
	\end{array}\right]
	\cdot
	\left[\begin{array}{llll}
		1 & 3000 & 33 & 5 \\
		1 & 6500 & 30 & 2 \\
		1 & 3500 & 31 & 4 \\
		1 & 2800 & 28 & 5
	\end{array}\right]\right)^{-1}
	\cdot
	\left[\begin{array}{llll}
		1 & 1 & 1 & 1 \\
		3000 & 6500 & 3500 & 2800 \\
		33 & 30 & 31 & 28 \\
		5 & 2 & 4 & 5
	\end{array}\right]
	\cdot
	\left[\begin{array}{c}
		4500 \\
		12000 \\
		8500 \\
		4000
	\end{array}\right].
	\end{small}
	$$

	Hence the coefficient for $\mathbf{X}$, i.e., $\Theta$, is
	$$
	\Theta = \left[\begin{array}{c}
		-22921 \\
		5 \\
		-210 \\
		2764
	\end{array}\right].
	$$

	\subsubsection{Standardization}
	We will give several definitions, which would help us standardize the whole process.
	\begin{defn}
		The cost function of multi-variable regression is:
		$$
		J\left(\theta_{0}, \theta_{1}, \cdots, \theta n\right)=\frac{1}{2 m} \sum_{m}^{i=1}\left[h_{\theta}\left(x^{i}\right)-y^{(i)}\right]^{2}
		$$
	\end{defn}
	\begin{defn}
		The predict function $h_\theta(x)$ is defined as:
		$$
		h_{\theta}(x)=\Theta^{T} X=\theta_{0} x_{0}+\theta_{1} x_{1}+\cdots+\theta_{n} x_{n}.
		$$
	\end{defn}
	The normal equation is to find the parameters that minimize the cost function by solving the following equation:
	$$
	\frac{{\partial}}{\partial \Theta_{j}} J\left(\theta_{j}\right)=0.
	$$
	\begin{defn}
		For $m$ training samples, and each sample has $n$ features, the dataset is defined as
		$$
		\mathbf{X}=\left[\begin{array}{ccc}
			x_{0}^{(1)} & \cdots & x_{n}^{(1)} \\
			\vdots & \ddots & \vdots \\
			x_{0}^{(\mathrm{m})} & \cdots & x_{n 0}^{(\mathrm{m})}
		\end{array}\right],
		$$
		where $x_{j}^{(i)}$ denotes the $j$-th feature of $i$-th sample.
	\end{defn}
	\begin{defn}
		The coefficient is defined as:
		$$
		\Theta=\left[\theta_{0}, \theta_{1}, \cdots, \theta_{n}\right]^{T}.
		$$
	\end{defn}
	\begin{defn}
		The output variable is defined as:
		$$
		\mathbf{Y}=\left[y^{(1)}, y^{(2)}, \cdots, y^{(m)}\right]^{T}.
		$$
	\end{defn}
	With the definitions above, we could define the cost function.
	\begin{defn}
		The cost function is defined as:
		$$
		\begin{aligned}
			J\left(\theta_{0}, \theta_{1}, \cdots, \theta_{n}\right)&=\frac{1}{2 m}(\mathbf{X} \cdot \Theta-\mathbf{Y})^{T}(\mathbf{X} \cdot \Theta-\mathbf{Y}) \\
			&=\frac{1}{2 m}\left(\mathbf{Y}^{T} \mathbf{Y}-\mathbf{Y}^{T} \mathbf{X }\Theta-\Theta^{T} \mathbf{X}^{T} Y+\Theta^{T} \mathbf{X}^{T} \mathbf{X} \Theta\right).
		\end{aligned}
		$$
	\end{defn}
	
	To find the derivative, we get:
	$$
	J'=\frac{1}{2 m}\left[\frac{\partial \left(\mathbf{Y}^{T} \mathbf{Y}\right)}{\partial \Theta}-\frac{\partial \left(\mathbf{Y}^{T} \mathbf{X} \Theta\right)}{\partial \Theta}-\frac{\partial \Theta^{T} \mathbf{X}^{T} \mathbf{Y}}{\partial \Theta}+\frac{\partial \Theta^{T} \mathbf{X}^{T} \mathbf{X} \theta}{\partial \Theta}\right].
	$$
	
	By simplifying we have:
	$$
	\frac{1}{2 m}\left(-2 \mathbf{X}^{T} \mathbf{Y}+2 \mathbf{X}^{T} \mathbf{X} \Theta\right)=0.
	$$
	
	Hence the coefficients could be derived as
	$$
	\Theta=\left(\mathbf{X}^{T} \mathbf{X}\right)^{-1} \mathbf{X}^{T} Y.
	$$
	
	We have ended our proof.
	
	% ---------------- end of document body---------------------
	%---------------- start of references------------------
	\clearpage
	\bibliographystyle{agsm} % this gives Harvard style references
	\bibliography{mybibfile} % change to the name of your .bib file

	%---------------- end of references------------------
\end{document}